**The Effects of Connectivity and Automation on Saturation Headway and Capacity at Signalized Intersections**

**Ali Hajbabaie**, Ph.D., Corresponding Author
Associate Professor
Civil, Construction, and Environmental Engineering Department
North Carolina State University
Fitts-Woolard Hall, 915 Partners Way, Raleigh, NC 27606
Tel: 919-515-5938; Email: ahajbab@ncsu.edu

**Mehrdad Tajalli**, Ph.D.
Graduate Research Assistant
Civil, Construction, and Environmental Engineering Department
North Carolina State University
Email: mtajal@ncsu.edu

**Eleni Bardaka**, Ph.D.
Assistant Professor, Department of Civil, Construction, and Environmental Engineering,
North Carolina State University
Fitts-Woolard Hall, 915 Partners Way, Raleigh, NC 27606
Tel: 919-515-8576; Email: ebardak@ncsu.edu

## ABSTRACT
This paper analyzes the potential effects of connected and automated vehicles on saturation headway and capacity at signalized intersections. A signalized intersection is created in Vissim as a testbed, where four vehicle types are modeled and tested: (I) human-driven vehicles (HVs), (II) connected vehicles (CVs), (III) automated vehicles (AVs), and (IV) connected automated vehicles (CAVs). Various scenarios are defined based on different market penetration rates of these four vehicle types. AVs are assumed to move more cautiously compared to human drivers. CVs and CAVs are supposed to receive information about the future state of traffic lights and adjust their speeds to avoid stopping at the intersection. As a result, their movements are expected to be smoother with a lower number of stops. The effects of these vehicle types in mixed traffic are investigated in terms of saturation headway, capacity, travel time, delay, and queue length in different lane groups of an intersection. A Python script code developed by Vissim is used to provide the communication between the signal controller and CVs and CAVs to adjust their speeds accordingly. The results show that increasing CV and CAV market penetration rate reduces saturation headway and consequently increases capacity at signalized intersections. On the other hand, increasing the AV market penetration rate deteriorates traffic operations. Results also indicate that the highest increase (80%) and decrease (20%) in lane group capacity are observed respectively in a traffic stream of 100% CAVs and 100% AVs.

## INTRODUCTION
Automated vehicles (AVs) are emerging with the promise of improving mobility and safety on highway facilities. Many research laboratories, vehicle manufacturers, and technology companies are currently researching and testing AVs on highway facilities. For instance, Waymo's AVs have traveled about 10 million miles on public roads in 25 cities across the US (1). Several car manufacturers such as Tesla, Cadillac, and Audi are building semi-automated commercial vehicles, while fully automated vehicles are expected to emerge by 2050 (2,3).

Vehicle connectivity is also expected to play an essential role in improving mobility and safety (4–7). Current studies show that AVs are programmed to behave conservatively, perhaps to reduce the likelihood of severe crashes in the absence of information from other vehicles and obstacles that are not visible to the sensors of AVs (8,9). Establishing dynamic communication among vehicles, infrastructure, and other wireless devices enables vehicles to collect real-time data and predict the future states of other users on the road more accurately. Consequently, this type of dynamic communication is expected to reduce the likelihood of crashes. As a result, connected AVs can drive more aggressively without increasing the risk of collision with other users (10). The advisory information also helps human drivers to oversight the upcoming traffic condition and adjust their speed and maneuvers appropriately (11–15).

This study aims to understand the potential effects of connectivity and automation on traffic operations at signalized intersections. Existing studies have paid attention to the operations of automated and connected vehicles on freeway facilities (16–20). Still, the interaction of automated and connected vehicles with other vehicles and the signal controller has not received the same amount of attention (21). We have considered four types of connectivity/automation to account for different driving behaviors and their interactions in traffic stream: I) human-driven vehicles (HVs), II) connected vehicles (CVs), III) automated vehicles (AVs), and IV) connected and automated vehicles (CAVs). Various scenarios are examined for different market penetration rates (proportion of each vehicle type in the traffic stream) of these vehicles.

The potential effects of CVs, AVs, and CAVs on intersection-level and lane-group level saturation headway, travel time, delay, and queue length are studied using a simulated testbed created in Vissim. Full automation is assumed for both AVs and CAVs. These four aformentioned vehicle types are simulated by changing car following model parameters in Vissim. Specifically, AVs receive information on their surrounding environment only from their on-board sensors and are assumed to operate with decision algorithms that are more conservative than HVs (9,22). CVs are assumed to receive information about the future state of traffic lights and adjust their speeds to avoid stopping at the intersection. As a result, the movements of CVs are expected to be smoother with a lower number of stops. CAVs are assumed to combine all the capabilities of AVs and CVs. A Python script code developed by Vissim is used to model the communication between the signal controller, CVs and CAVs to adjust their speed accordingly.

In the remainder of the paper, we first summarize the related literature and introduce the methodology. Then, the case study and numerical results are presented. The discussion of results is followed by concluding remarks and trends for future research.

## LITERATURE REVIEW

The driving behavior associated with HVs, CVs, AVs, and CAVs differs by automation and connectivity level. Lower automation levels aim to assist human drivers through technologies enabled by onboard computers and sensors, such as adaptive cruise control (23), collision warning (24,25), collision avoidance (26), or assistant braking (27). On the other hand, higher automation levels enable AVs and CAVs to take complete control of the vehicle's movements without any assistance from the human driver by predicting the future trajectory of surrounding vehicles and avoiding any potential collisions.

The interaction between vehicles with different levels of connectivity and automation will be a challenge in the near future since these vehicles have different driving behaviors (18). AVs, as of now, are programmed to behave cautiously while interacting with HVs (28,29). Human drivers require a higher reaction time to respond to any changes in the driving environment. Therefore, AVs need to consider various decision scenarios to overcome the uncertainty associated with human driver decisions (30). Sadigh et al. (*16*) showed that AVs could influence human driver behavior and yield a more efficient performance. Sezer et al. (31) also clarified that making AVs more aggressive could yield to operating higher traffic volumes in a mixed autonomy environment without compromising the safety when there is communication between vehicles.

Connectivity can further improve the efficiency and reliability of automated systems (as well as human-controlled systems). The information sharing between vehicles and infrastructure increases the chance of reliable decisions during driving, especially with respect to car following and lane changing (32). In addition, CAVs can improve traffic mobility without sacrificing safety. For instance, controlling the trajectory of CAVs upstream of signalized intersections based on advanced knowledge of signal phase and timing (SPaT) increases intersection throughput and reduces the experienced delay and risk of collisions among vehicles (25,33–38). Moreover, the trajectory of CAVs can be managed to avoid stops at the intersection and minimize fuel consumption (12,34,39).

While many studies have examined the possible effects of connected and automated vehicles on traffic operations on uninterrupted flow facilities (16–20), the impacts of connectivity and automation on interrupted flow facilities, especially signalized intersections, are not thoroughly studied. Existing studies focus on either using signal timing information to plan the arrival of CAVs (33–35), jointly optimizing signal timing parameters and CAV trajectories (11,40–44), or designing a signal-free environment with a fleet of 100% CAVs (45–55). The effects of different market penetration levels of connectivity and automation on the saturation headway and capacity at signalized intersections are unknown. This study aims to fill this gap and provide insights into how a different mix of CVs, AVs, and CAVs in the traffic stream will influence saturation headway and capacity on signalized intersections.

## METHODOLOGY

In this study, various driving behaviors are defined in Vissim simulation environment to represent HVs, CVs, AVs, and CAVs. Vissim is calibrated for the base condition representing a fleet of 100% HVs on an exclusive through lane group. Different combinations of HV, CV, AV, and CAV market penetration rates are created, and the saturation headway along with delay, travel time, queue length, and throughput are measured for each scenario. The saturation headway values are used to fit a model of saturation headway as a function of HV, CV, AV, and CAV market penetration rate, lane group configuration, and turning percentage. The outcome of this model is used to determine capacity adjustment factors. Additional information is provided in the following sections.

### Vissim Simulation

A microscopic simulation testbed is developed in Vissim to study the connectivity and automation technologies' effects on traffic operations at signalized intersections. The simulation testbed provides the ability to consider various driving behaviors associated with connectivity and automation. Vissim allows simulating connected and automated vehicles and their interaction with conventional human-driven vehicles. In addition, information exchange between vehicles and the infrastructure is modeled through Vissim's Component Object Model (COM) interface. Finally, Vissim provides a host of outputs ranging from vehicle-level output to network-level performance measures.

### Driving Behavior

The driving behavior of vehicles with full automation is adopted from existing practical studies on the behavior of AVs. In particular, the findings of the CoEXist project are used to simulate the movement of AVs in Vissim. The recommendations of the CoEXist project are based on the empirical analysis of data collected in the Netherlands. The experimental results are confirmed by Vedecom Tech and several simulation tests done by PTV Group (17).

#### *Car-Following Behavior*

The CoEXist project recommended three variants of driving models: 1) CoEXist cautious model, 2) CoEXist normal model, and 3) CoEXist all-knowing model (17). The cautious driving behavior respects the road code and always ensures moving safely on the road. There is always a brick wall distance between a cautious-driving vehicle and its immediate leading car. This means that if the leading vehicle comes to an stantanous stop, the cautious -driving vehicle can stop as well and avoided a crash. In addition, a large gap is required to perform a lane change maneuver or pass an un-signalized intersection. The normal driving behavior is very similar to the behavior of a human driver with the additional capacity of measuring distances and speeds of surrounding vehicles by collecting information from sensors. All-knowing driving behavior assumes a perfect perception of the surrounding environment and receives vehicle-to-vehicle and vehicle-to-infrastructure communications. This driving behavior is associated with smaller gaps for all maneuvers. In this study, the movement of AVs is assumed to follow the CoEXist cautious model while the movement of CAVs is assumed to follow the CoEXist all-knowing model. The movement of HVs and CVs follow normal driving behavior with a difference in following variation, which is smaller for CVs. The signal timing information is shared with connected human-driven vehicles and connected automated vehicles. As a result, the driving behavior of vehicles under these two types will be different than the cases that the information is not received, as follows in the next section. TABLE 1 summarizes the calibrated components of Wiedemann 99 car-following parameters adopted from the CoEXist project.

**TABLE 1** Wiedemann 99 car-following model Calibration Components (CC) –Reference: (17)

| **Parameters** | **Definition** | **HV** | **CV** | **AV** | **CAV** |
|---|---|---|---|---|---|
| CC0 | Stand still distance ($m$) | 1.5 | 1.5 | 1.5 | 1 |
| CC1 | Headway time ($s$) | 1.6 | 1.6 | 2.2 | 1 |
| CC2 | Following variation ($m$) | 4 | 2 | 0 | 0 |
| CC3 | Threshold for entering following ($s$) | -8 | -8 | -10 | -6 |
| CC4 | Negative following threshold ($m/s$) | -0.35 | -0.35 | -0.1 | -0.1 |
| CC5 | Positive following threshold ($m/s$) | 0.35 | 0.35 | 0.1 | 0.1 |
| CC6 | Speed dependency of oscillation ($1/(m/s)$) | 11.44 | 11.44 | 0 | 0 |
| CC7 | Oscillation acceleration ($m/s^2$) | 0.25 | 0.25 | 0.1 | 0.1 |
| CC8 | Standstill acceleration ($m/s^2$) | 3.5 | 3.5 | 2 | 4 |
| CC9 | Acceleration with 50 mph ($m/s^2$) | 1.5 | 1.5 | 1.2 | 2 |

*Signal Control Behavior*

When there is no information available about the future signal timing plan at an intersection, vehicles either follow their lead vehicles or travel at their own desired speed. If they hit the green signal, they will go through the intersection; otherwise, they stop for the red light. Sharing signal timing information with upcoming vehicles can change their driving behavior as they approach the intersection. TABLE 2 shows Vissim's car following behavior for different connectivity and automation levels as vehicles arrive at an intersection. HVs and AVs need to continuously check the signal timing status during the yellow time to avoid red-light violations. Although CVs receive the signal timing plans, the human driver still needs to check to constantly ensure safe entrance to the intersection. On the other hand, CAVs receive information on future signal timing plans and do not need to check them continuously. In addition, the safety factor for AVs and CAVs is considered higher to ensure no collision will occur with other vehicles at the intersection.

**TABLE 2** Car-following behavior near signalized intersections–Reference: (17)

| Attribute | HV | CV | AV | CAV |
|---|---|---|---|---|
| behavior at the amber signal | continuous check | continuous check | continuous check | one decision |
| behavior at red/amber signal | go | go | stop | stop |
| reduced safety distance factor | 0.6 | 0.6 | 1 | 1 |
| reduced safety start upstream of stop line (m) | 100 | 100 | 100 | 100 |
| reduced safety end upstream of stop line (m) | 100 | 100 | 100 | 100 |

*Other Behaviors*

Other driving behavior parameters suggested by the CoEXist project are shown in TABLE 3. Enforce absolute braking distance (EABK) is active for AVs since they drive cautiously on the road. Based on EABK, a further gap between the following and leading vehicles is maintained to allow AVs to stop safely anytime, even if the lead vehicle stops instantly. Vissim does not consider any stochasticities associated with automated vehicles. CAVs can interact with more than one vehicle in the traffic stream, but non-CAVs interact with the most immediate vehicle.

**TABLE 3** Other driving behaviors for various automation levels–Reference: (17)

| Driving logic | Enforce absolute braking distance (EABK) | Use implicit stochastics | Number of interaction vehicles | Increased desired acceleration |
|---|---|---|---|---|
| HV | OFF | ON | 1 | 100-110% |
| CV | OFF | ON | 1 | 100% |
| AV | ON | OFF | 1 | 100% |
| CAV | OFF | OFF | >1 | 110% |

## Vissim Calibration

Wiedemann's 99 car-following model is selected because it can model both automated and connected automated vehicles (17). The main reason for using this model is that the model parameters for AVs and CAVs are found for this model and not Wiedemann's 74 car-following model. Model parameters are calibrated to match the saturation headway. The headway depends on two main factors in Vissim: 1) the desired speed and 2) the car-following characteristics. The desired speed is defined as the speed that vehicles use in free flow conditions. Since the desired speed is assumed fixed for all vehicle types, only the car-following parameters should be calibrated. As shown in TABLE 1, the Wiedemann 99 car-following model contains ten parameters. However, only two parameters (i.e., CC0 and CC1) influence the intersection headway significantly (56). CC0 is the average desired distance between two vehicles in meters at a stand-still while queuing before the traffic signal. The headway (CC1) describes the speed-dependent part of the safety distance a driver desires. Therefore, various combinations of these two factors are tested to calibrate the model.

## Advisory Speeds

CVs and CAVs can adjust their speeds based on the received information on future signal timing plans to smoothen their movement and arrive at the intersection during the green signal. This research utilizes a Python script code developed by the PTV group to allow communications between the signal controller and CVs and CAVs. The script adjusts the minimum and maximum speed required to arrive at the intersection during a green light. If the minimum speed is less than the desired speed, the vehicle moves with the desired speed; otherwise, a constant smooth speed will be provided. It should be noted that the constant speed needs to be higher than a certain amount to avoid crawling. Vissim's default value is five mph, which is used in this study. Similarly, the maximum speed to the intersection is calculated and compared to the desired speed. If the maximum speed is greater than the desired speed, the desired speed will be used. Otherwise, the maximum speed is considered as the optimal speed for arriving at the intersection during the green signal.

**Saturation Headway and Capacity Adjustment Factors**

The saturation headway is defined as the average headway between the fourth and tenth passenger cars in the queue when the traffic light changes from red to green. The headways are quantified using a script that analyzes the trajectory of vehicles and determines the exact time their front bumpers hit the stop line. This approach is based on Chapter 6 of the ITE Manual of Transportation Studies (MTES) (57). The estimated saturation flow rate from field data is defined as the number of vehicles that pass the stop bar at the intersection with no interruptions in a defined period. Based on the MTES approach for collecting saturation flow rate data, a timer is started when the fourth vehicle in a queue passes the stop bar; at that point, the queue typically begins maintaining consistent headway after any incurred start-up lost time. The timer is stopped when either the seventh, eighth, ninth, or tenth vehicle passes the stop bar, whichever is the last in the stopped queue. For a typical standard deviation in saturation flow rate of 140 vehicles per hour, the MTES suggests observing a minimum of 30 valid queues to estimate the mean saturation flow rate within 50 vehicles per hour of the true rate and a 95% confidence level (57). We note that any queue with less than seven passenger cars would be excluded from this analysis to avoid confounding effects. Saturation headway data is aggregated at the 15-minute period to be consistent with the 15-minute analysis period of the highway capacity manual (HCM) (58).

A regression approach can use either additive or multiplicative models to predict saturation headway. This study employes additive models due to their simplicity and consistency with the HCM. In additive models, the base saturation headway (all HVs on the exclusive through lane group) is adjusted by subtracting or adding values corresponding to different market penetration rates of CVs, AVs, and CAVs. In addition, the saturation headway is adjusted for exclusive right or left-turn lanes. The intercept ($h_s$) represents base saturation headway, and the dependent variable ($h_s^{adj}$) represents adjusted saturation headway under a combination of CVs, AVs, and CAVs.

$$h_S^{adj} = h_s + \beta_1(\mathrm{CV}) + \beta_2(\mathrm{AV}) + \beta_3(\mathrm{CAV}) + \beta_4(\mathrm{D_{EXL}}) + \beta_5(\mathrm{D_{EXR}}) + \beta_6(\mathrm{D_{SHTR}}) + \beta_7(\mathrm{D_{SHTR}})(\mathrm{RT}) \quad (1)$$

where:

| | |
|---|---|
| $h_s^{adj}$: | adjusted saturation headway (seconds), |
| $h_s$ : | base saturation headway for non-work zone conditions, |
| CV : | market penetration rate of CVs in real number, |
| AV : | market penetration rate of AVs in real number, |
| CAV : | market penetration rate of CAVs in real number, |
| $\mathrm{D_{EXL}}$ : | Exclusive left turn lane (0: no, 1: yes), |
| $\mathrm{D_{EXR}}$ : | Exclusive right turn lane (0: no, 1: yes), |
| $\mathrm{D_{SHTR}}$ : | Shared through and right lane group (0: no, 1: yes), and |
| RT : | Right turn percentage on shared lanes (%). |

The capacity of each lane group is a function of the saturation headway, effective green, and cycle length as follows:

$$c = \left(\frac{3600}{h_s}\right)\left(\frac{g}{C}\right) \tag{2}$$

where:

| | |
|---|---|
| $c$ : | capacity (passenger cars per hour per lane), |
| $g$ : | effective green (seconds), and |
| $C$ : | cycle length (seconds). |

Therefore, the adjusted capacity $c^{adj}$ can be determined using the following equation:

$$c^{adj} = \left(\frac{3600}{h_s^{adj}}\right)\left(\frac{g}{C}\right) \tag{3}$$

Assuming the signal timing parameters are unchanged due to the presence of CVs, AVs, and CAVs, capacity adjustment factors (CAF) can be determined as:

$$\text{CAF} = \frac{c^{adj}}{c} = \frac{h_s}{h_s^{adj}} \tag{4}$$

Note that capacity adjustment factors can be larger than one as a fleet of 100% CAVs is expected to increase the capacity.

## INTERSECTION TESTBED

FIGURE 1 shows the layout of the intersection testbed used in this research. The testbed includes various lane groups. The eastbound approach has exclusive left turn, through, and right turn lane groups. Other approaches include a shared right turn and through lane group. Fixed-time signal timing is used, and the signal timing parameters are optimized using Vistro (59). The demand for the eastbound entry is 900 vehicle/hour, and the demand for other approaches is 1200 vehicle/hour. The turning percentage for left-turn movement is 15% for all approaches. The right-turn percentage for eastbound, northbound, westbound, and southbound are assumed to be 15%, 5%, 15%, and 25%, respectively. Six different market penetration rates for CVs, AVs, and CAVs are simulated: 0%, 20%, 40%, 60%, 80%, and 100%. For instance, a scenario can have 20% HVs, 20% CVs, 40% AVs, and 20% CAVs. In total, 56 scenarios are considered. Each scenario is run ten times with different random seeds to account for randomness.

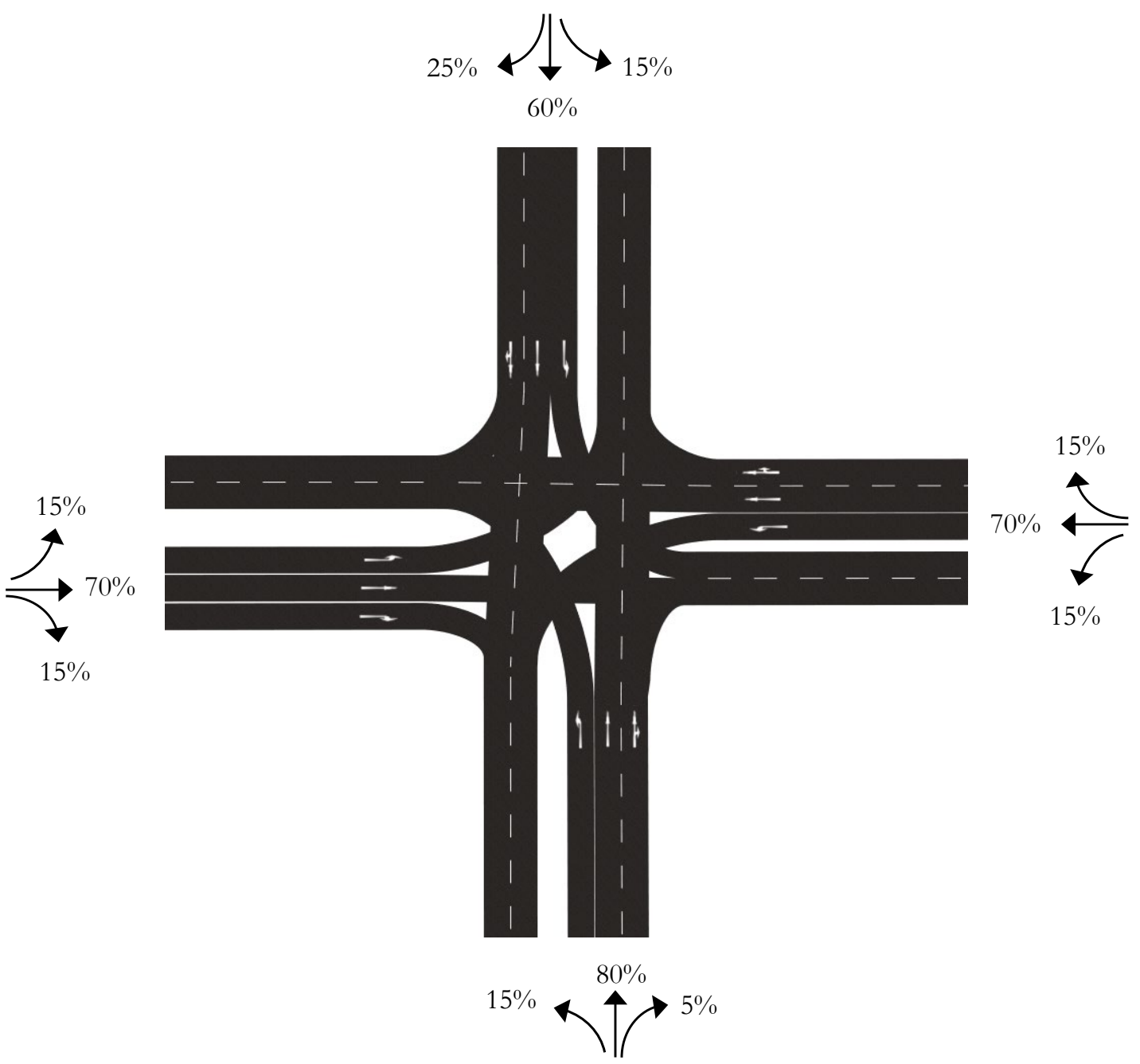

**FIGURE 1** Case study intersection

## RESULTS

### Lane Group-Level Analysis

*Saturation Headway*

FIGURE *2* shows the saturation headway of CVs, AVs, and CAVs compared to the base case, i.e., 100% HVs (shown as 0% penetration rate in the figure). The results are shown for exclusive right turn, through, and left turn lane groups. Increasing the market penetration rate of CVs and CAVs decreases the saturation headway at all exclusive lanes. This trend could be attributed to having access to advanced information about the future signal plan of the traffic light. As a result, the drivers of CVs are ready to start to move with shorter start-up lost time and reaction times. CAVs have shorter saturation headway than CVs due to the lower levels of uncertainty associated with the absence of human driving.

In contrast to CVs and CAVs, a higher penetration rate of AVs is associated with an increase in the saturation headway. This trend is primarily because AVs are programmed to travel more cautiously near the intersection to avoid collisions. As expected, the saturation headway is lower on exclusive through lanes compared to the exclusive left- and right-turn lanes because vehicles need to slow down to negotiate the curve.

A similar analysis is performed for shared right-turn and through lanes. The trends are similar to exclusive right-turn lanes. However, different right turn percentages are not found to substantially impact the saturation headway in the shared right and through lanes.

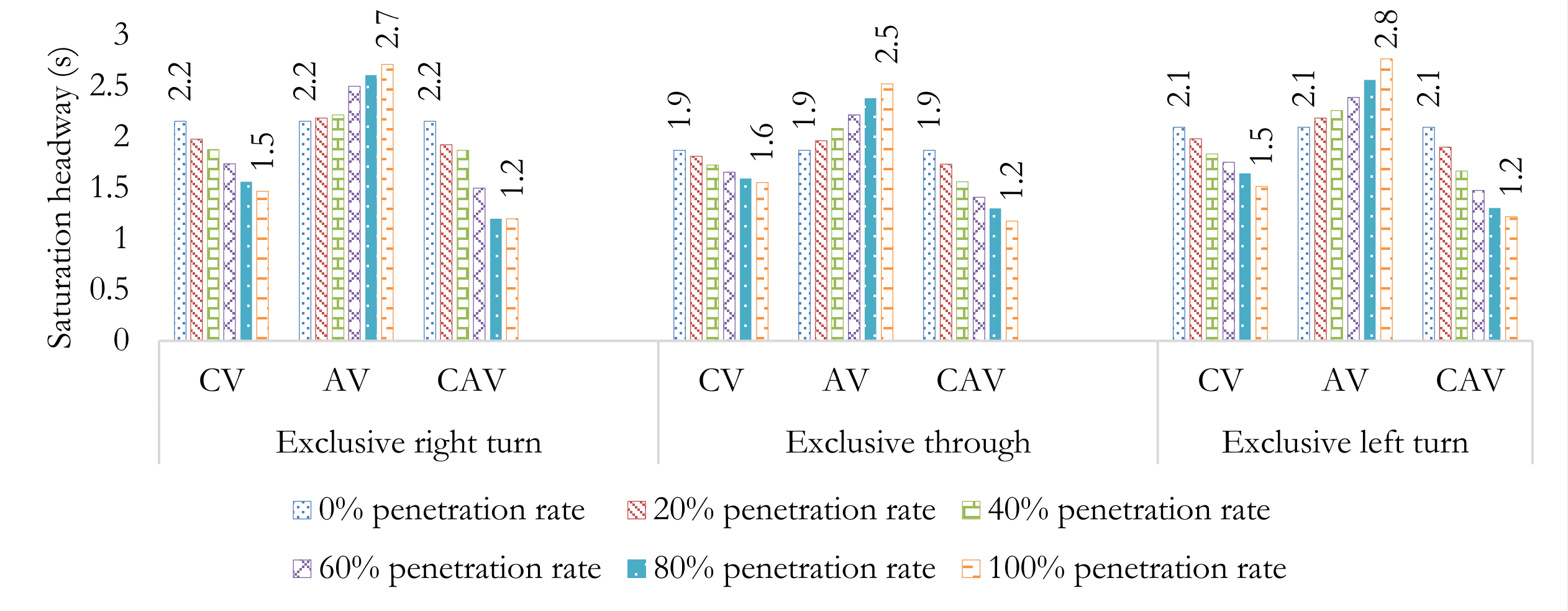

**FIGURE 2** Saturation headway for exclusive right, through, and left turning lanes

Multi-linear regression analysis is performed to predict the saturation headway as a function of CV, AV, and CAV market penetration and lane configuration. The base condition in this model is having an exclusive through lane with 100% HVs. Indicator variables specify different lane configurations corresponding to exclusive left, exclusive right, and shared through and right turn lanes. The shared through and right turn lane coefficient was close to zero and had a large p-value (>0.85). Therefore, it was removed from the model. Table 4 presents the final model results, and Eq. (5) presents the estimated model.

**TABLE 4** Regression Results

| | *Coefficients* | *Standard Error* | *t Stat* | *P-value* | *Lower 95%* | *Upper 95%* |
|---|---|---|---|---|---|---|
| Intercept ($h_s$) | **1.95** | 0.02 | 95.58 | 6.1381E-180 | 1.91 | 1.99 |
| CV | **-0.51** | 0.03 | -17.19 | 1.92367E-42 | -0.57 | -0.45 |
| AV | **0.56** | 0.03 | 18.69 | 3.58368E-47 | 0.50 | 0.62 |
| CAV | **-0.91** | 0.03 | -30.47 | 1.55072E-80 | -0.97 | -0.85 |
| Exclusive Left ($D_{EXL}$) | **0.11** | 0.02 | 6.86 | 6.94443E-11 | 0.08 | 0.14 |
| Exclusive Right ($D_{EXR}$) | **0.11** | 0.02 | 6.86 | 7.1267E-11 | 0.08 | 0.14 |

$$h_s^{adj} = 1.95 - 0.51(\text{CV}) + 0.56(\text{AV}) - 0.91(\text{CAV}) + 0.11\,(D_{EXL} + D_{EXR}) \qquad (5)$$

The model is associated with an adjusted R-squared value of 0.926 and passes the multi-linear regression assumptions. All the estimated model parameters are intuitive and have reasonable values. The intercept of the model in the base condition is 1.95, which means that the saturation headway in an exclusive through lane with 100% HVs is 1.95 seconds. The saturation headway decreases with an increase in the CV market penetration rate. Each additional percentage of CVs reduces the saturation headway by 0.0051 seconds. The opposite is found for AVs, where each additional percentage of them increases saturation headway by 0.0056 seconds. As expected, CAVs reduce the saturation headway. The model also suggests that the saturation headway in exclusive right or left turn lanes is 0.11 seconds longer than in an exclusive through lane regardless of the CV, AV, and CAV market penetration rate. The model helps estimate the saturation headway for a combination of CV, AV, and CAV market penetration rates. For instance, the saturation headway on an exclusive through lane with 10% HVs, 15% CVs, 25% AVs, and 50% CAVs can be calculated as $1.95 - 0.51(0.15) + 0.56(0.25) - 0.91(0.5) = 1.59$ seconds. While the model has a high adjusted R-squared value and the parameters are reasonable and according to our initial expectations, the output values of the model should be used cautiously. The main reason is that this study implements assumptions and changes, informed by the literature, in certain parameters of Vissim's car-following and lane-changing models, which were originally designed to represent human driving behavior. The primary purpose of the model is to illustrate trends.

The fitted model is used to visualize the effects of CV, AV, and CAV market penetration rates on saturation headway in exclusive through lane groups. **FIGURE 3** shows four heatmaps for 0%, 20%, 40%, and 60% HV market penetration rates in parts a through d, respectively. The x-axis shows the CAV market share, the y-axis shows the AV market share, and the CV market share is simply found by $\text{CV} = \text{HV} - \text{AV} - \text{CAV}$, where HV represents the human-driven vehicle market share. The trends are as expected: an increase in the CAV market share decreases the saturation headway, while the opposite trend is observed for AVs. Note that only in part (a), where the HV market share is zero, the market share of CVs, AVs, and CAVs can reach 100%.

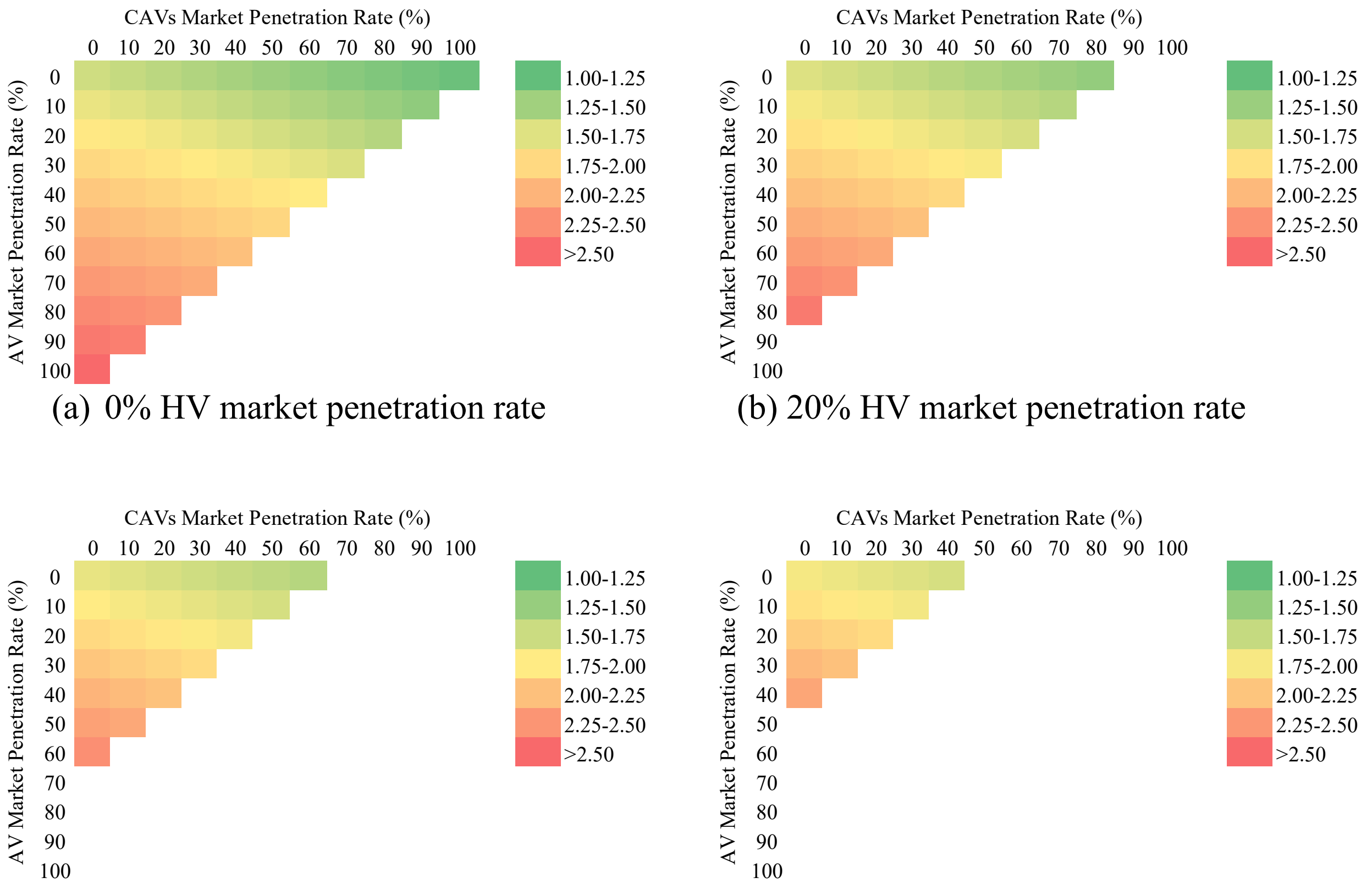

(a) 0% HV market penetration rate

(b) 20% HV market penetration rate

(c) 40% HV market penetration rate

(d) 60% HV market penetration rate

**FIGURE 3** Saturation headway (in seconds) for different CV, AV, and CAV market penetration rates on exclusive through lanes. [Note: $\text{CV} = 1 - \text{HV} - \text{AV} - \text{CAV}$]

*Capacity Adjustment Factors*

The capacity adjustment factor can be found using equation (6) shown below. This equation is directly found by dividing the base saturation flow rate by the adjusted saturation flow rate when CVs, AVs, or CAVs are present in the traffic stream. Note that any combination of market penetration rates can be entered into the equation as long as the market penetration rates of HVs, CVs, AVs, and CAVs add up to one.

$$CAF = \frac{h_s}{h_s^{adj}} = 1.95\big(1.95 - 0.51(\mathrm{CV}) + 0.56(\mathrm{AV}) - 0.91(\mathrm{CAV}) + 0.11\ (\mathrm{D_{EXL}} + \mathrm{D_{EXR}})\big)^{-1} \quad (6)$$

Similar to the previous section, the equation is used to visualize the effects of CV, AV, and CAV market penetration rates on capacity adjustment factors (see **FIGURE 4**). Note that it is assumed that the signal timing parameters are unchanged in the presence of connected and automated vehicles, and the change in capacity is only due to changes in saturation headway. Increasing the market penetration rate of CVs and CAVs increases the capacity of lane groups. On the contrary, more AVs decrease the capacity and bring it below the base capacity when only HVs are present in the traffic stream.

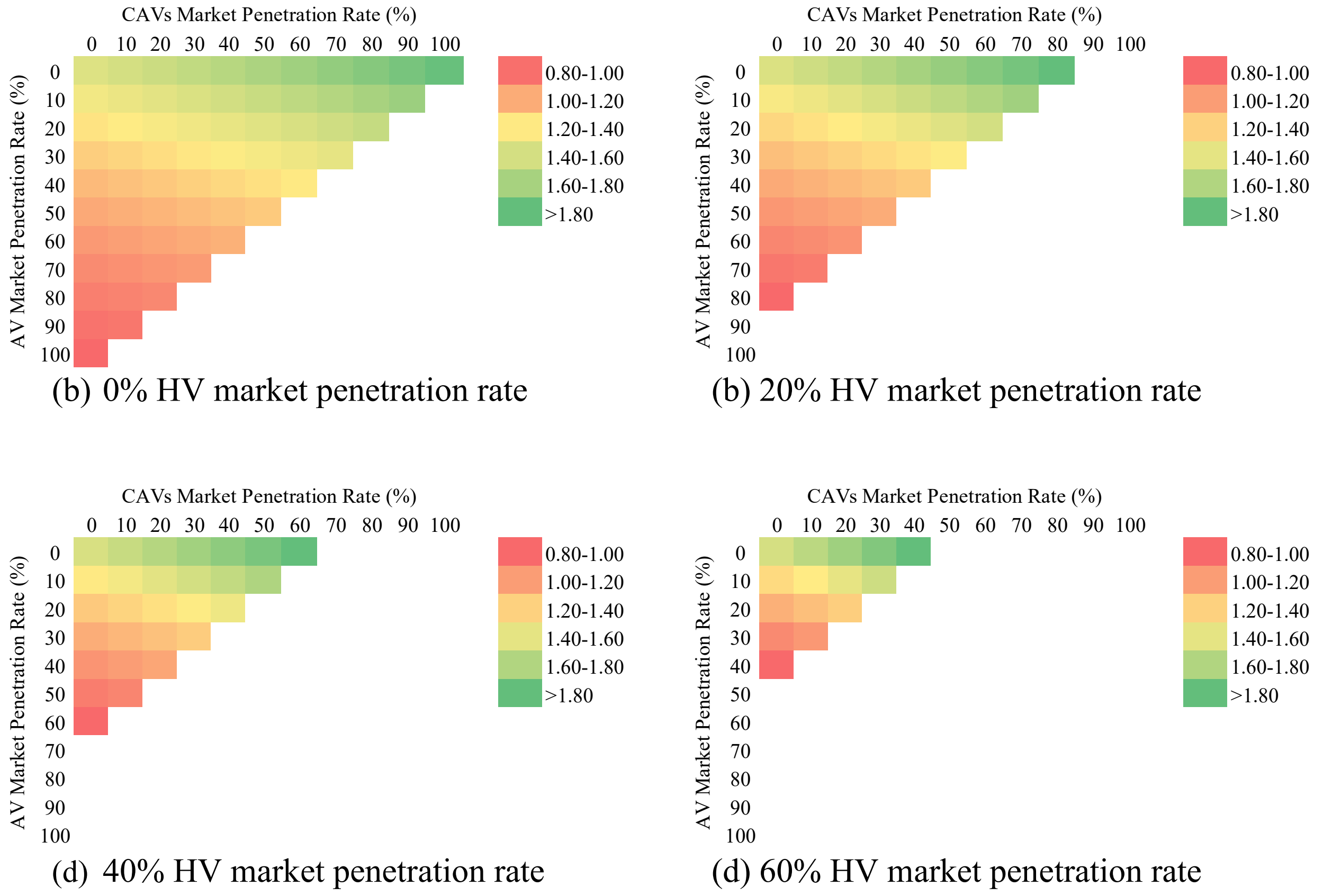

**FIGURE 4** Capacity Adjustment Factor for different CV, AV, and CAV market penetration rates on exclusive through lanes. [Note: $\mathrm{CV} = 1 - \mathrm{HV} - \mathrm{AV} - \mathrm{CAV}$]

*Average Delay*

In addition to saturation headway, the average delay of vehicles at each lane group under different market penetration rates and lane configuration is determined. FIGURE 5 shows that increasing the market penetration rate of CVs and CAVs leads to decreasing delay. On the other hand, increasing the AV market share increases the delay due to AVs' cautious driving behavior in the vicinity of the intersection. A significant increase in average delay on the exclusive through lane group is observed when the market penetration rate of AVs goes from 80% to 100%. The main reason for this steep increase is that the capacity of this lane group falls below the incoming volume when the AV market penetration rate reaches 100%.

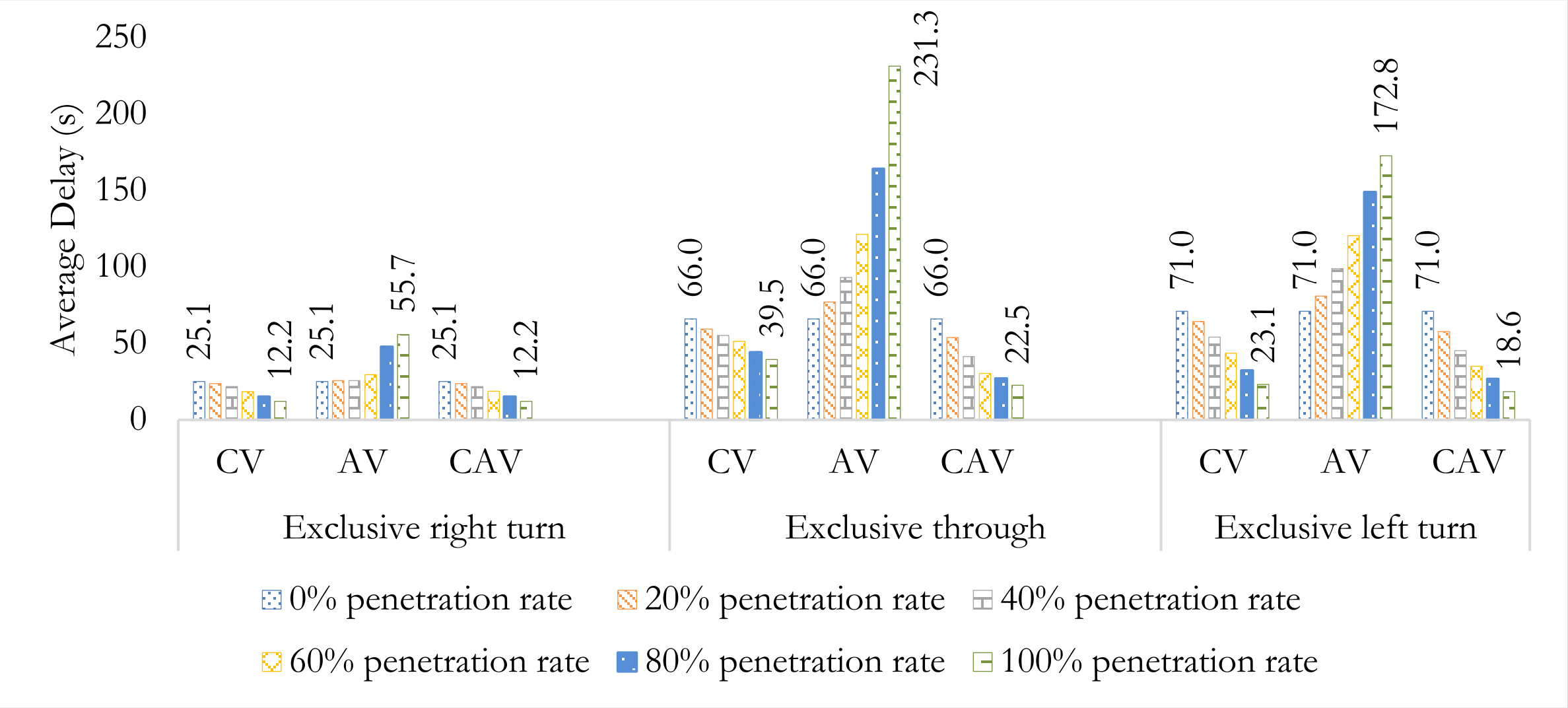

**FIGURE 5** Average delay for exclusive right, through, and left turning lane groups

FIGURE 6 shows the average delay of vehicles in shared lanes with different turning percentages. Similar trends are observed: increasing the market penetration rates of CVs and CAVs decreases the average delay, while increasing the AV market share is associated with an increase in the average delay. These trends are as expected and directly due to the trends that were observed before with regard to the effects of CVs, AVs, and CAVs on saturation headway and capacity.

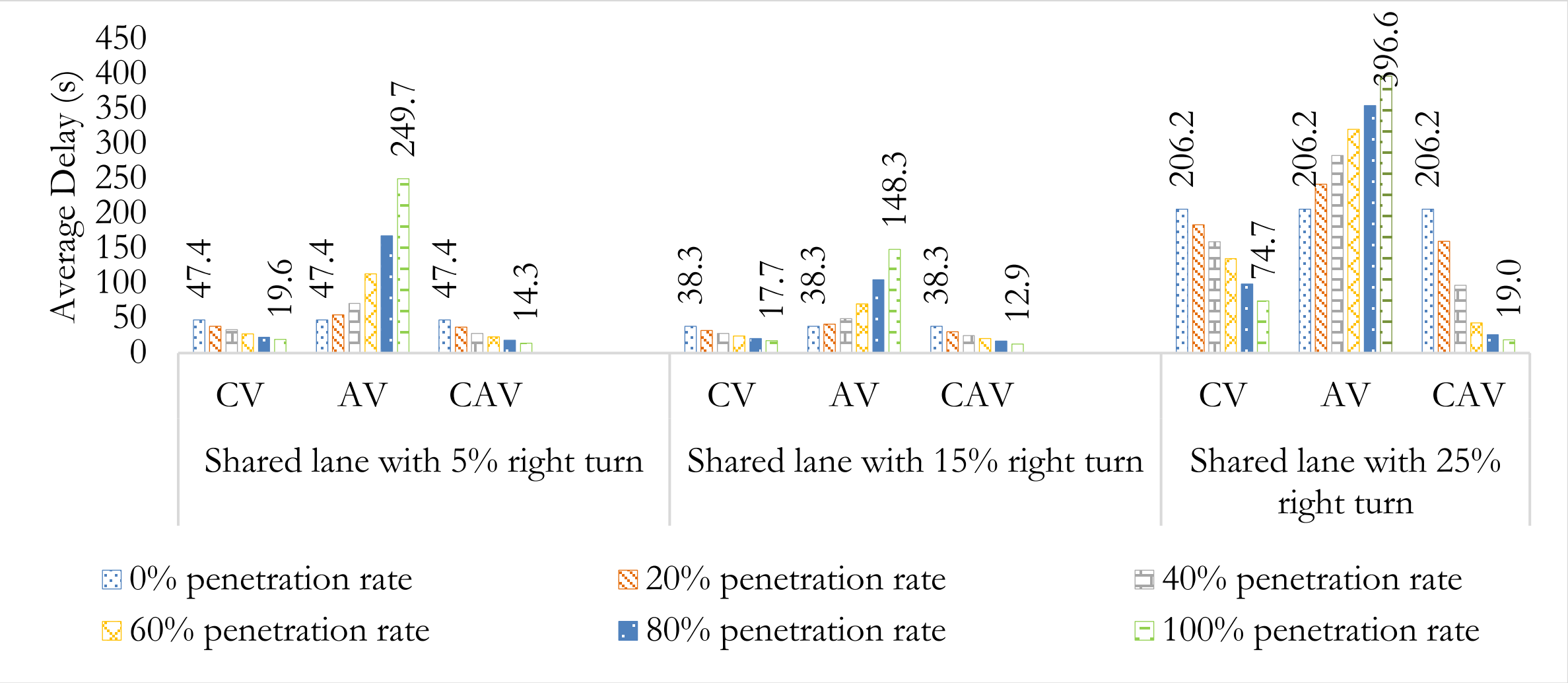

**FIGURE 6** Average delay for shared right and through lanes with 5%, 15%, and 25% turning percentages

*Queue Length*

FIGURE 7 shows the average queue length for each lane group of the intersection, separately. Increasing the market penetration rate of CVs and CAVs decreases the average queue length for all lane groups. The reduction rate is more substantial for through movements since they have higher demand volumes. On the other hand, increasing the market penetration rate of AVs results in a higher queue length, primarily due to the longer headways that AVs maintain compared to other vehicle types. We also observe that the queue length for southbound through (SBT) and southbound right (SBR) lane groups are longer than others due to high right-turning percentages (i.e., 25% right turn).

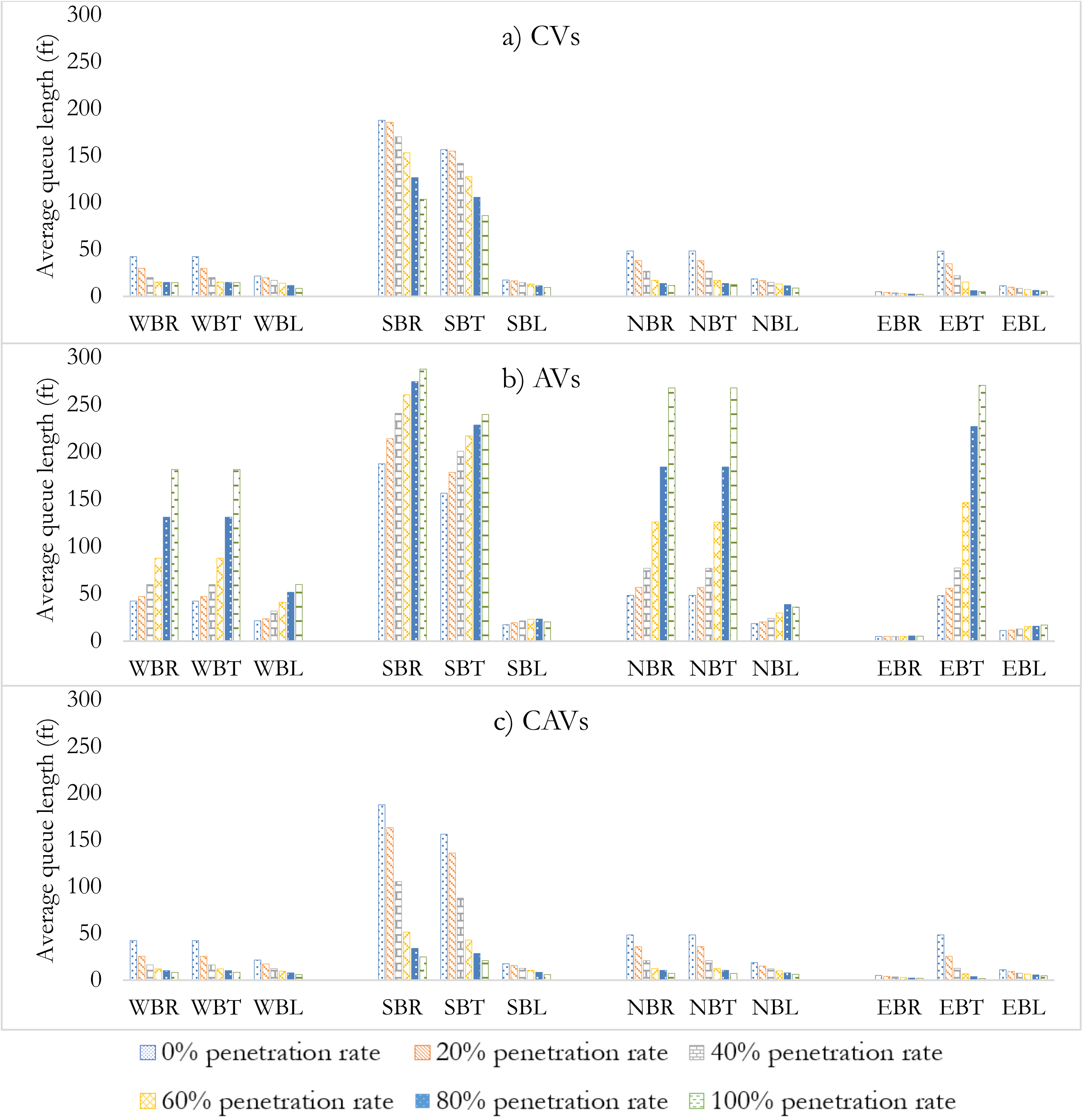

**FIGURE 7 Average queue length for all movements of the intersection**

### Intersection-Level Analysis

In addition to the lane group-level analysis, we analyze the effects of connectivity and automation on mobility performance measures at the intersection level. FIGURE 8 shows the average delay for the entire intersection. We find that increasing the penetration rate of CVs and CAVs decreases the average delay. However, increasing the AV penetration rate increases the average delay at the intersection. We observe that the lowest delay is associated with the highest number of CAVs in the intersection. Similar trends are also observed for the intersection average travel time, which is shown in FIGURE 9.

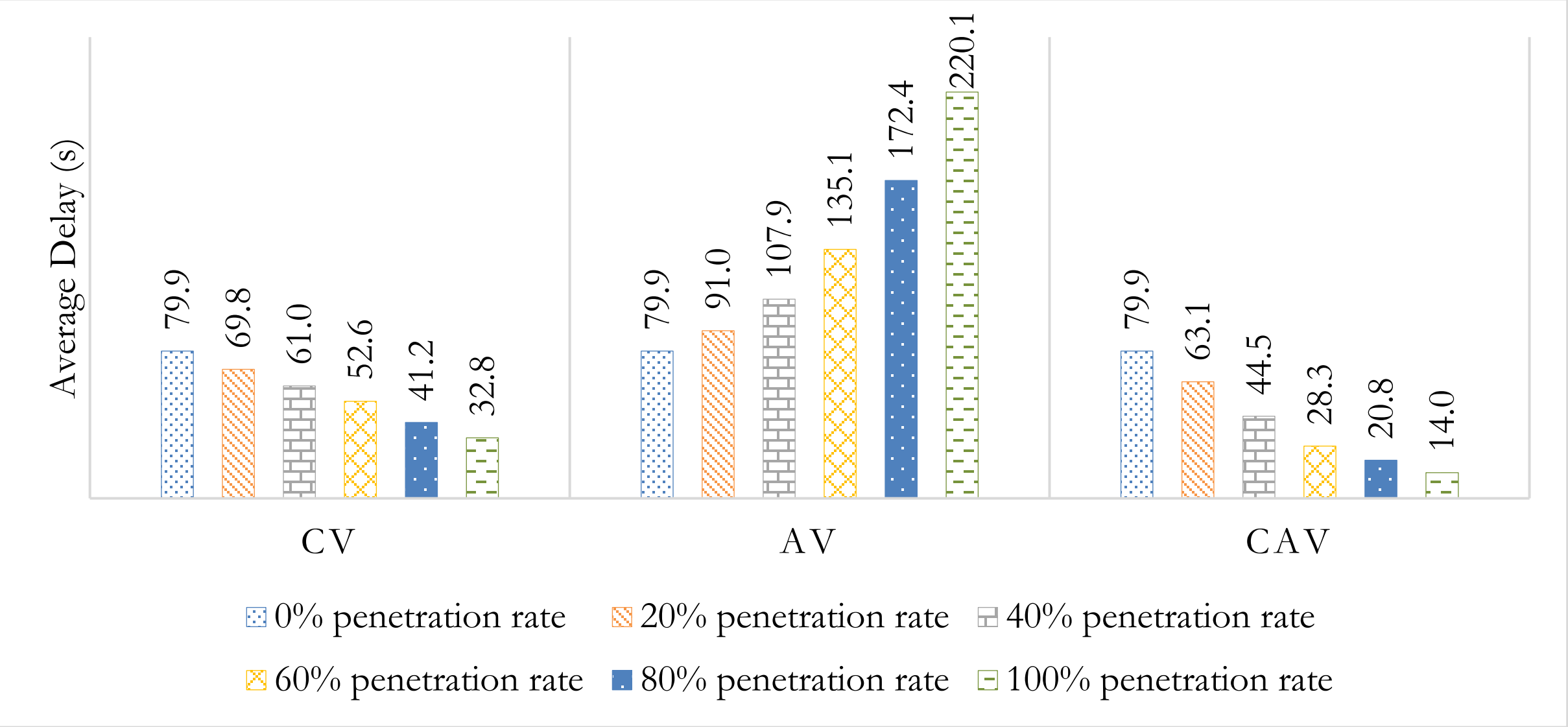

**FIGURE 8 Intersection-level average delay**

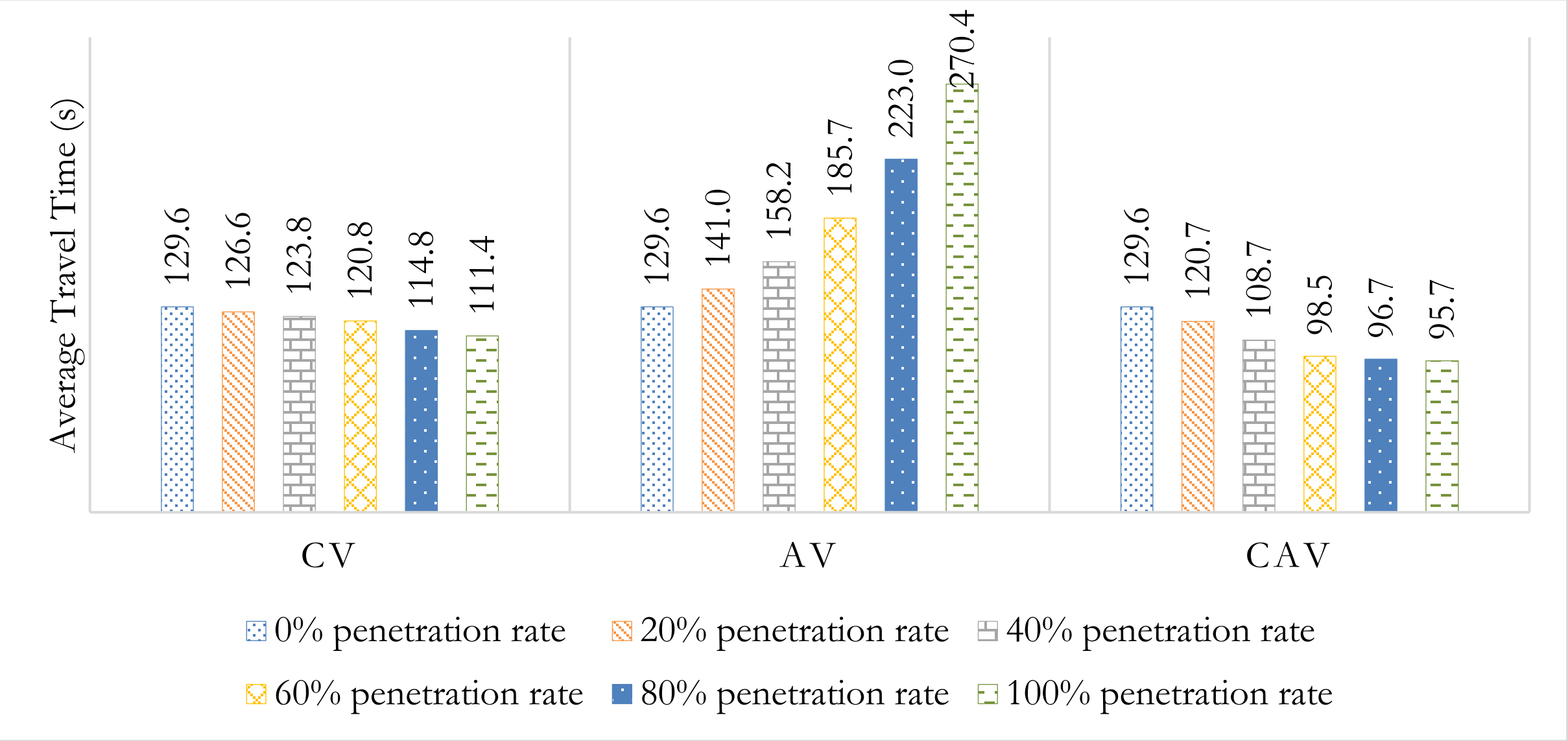

**FIGURE 9 Intersection-level average travel time**

Based on the results shown in FIGURE 10, increasing the penetration rate of CVs and CAVs increases the intersection throughput slightly. However, increasing the AVs penetration rate decreases intersection throughput substantially, which can be due to a reduction in intersection capacity. The maximum throughput for one hour of the simulation was 4,375 vehicles with 100%

CAVs in the traffic stream. On the other hand, the lowest throughput was 3,763 vehicles associated with 100% AV in the traffic stream.

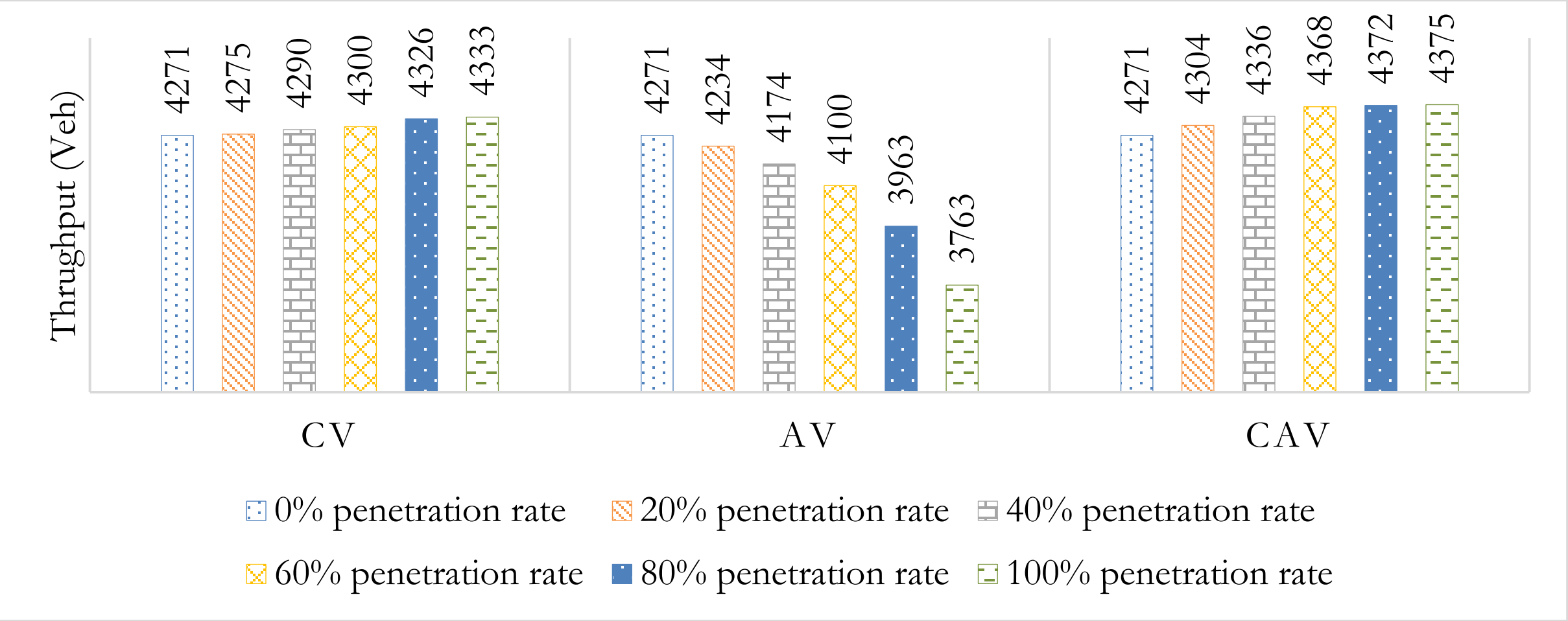

**FIGURE 10** Intersection throughput

FIGURE 11shows heatmaps for intersection-level saturation headway for different CV, AV, and CAV market shares.

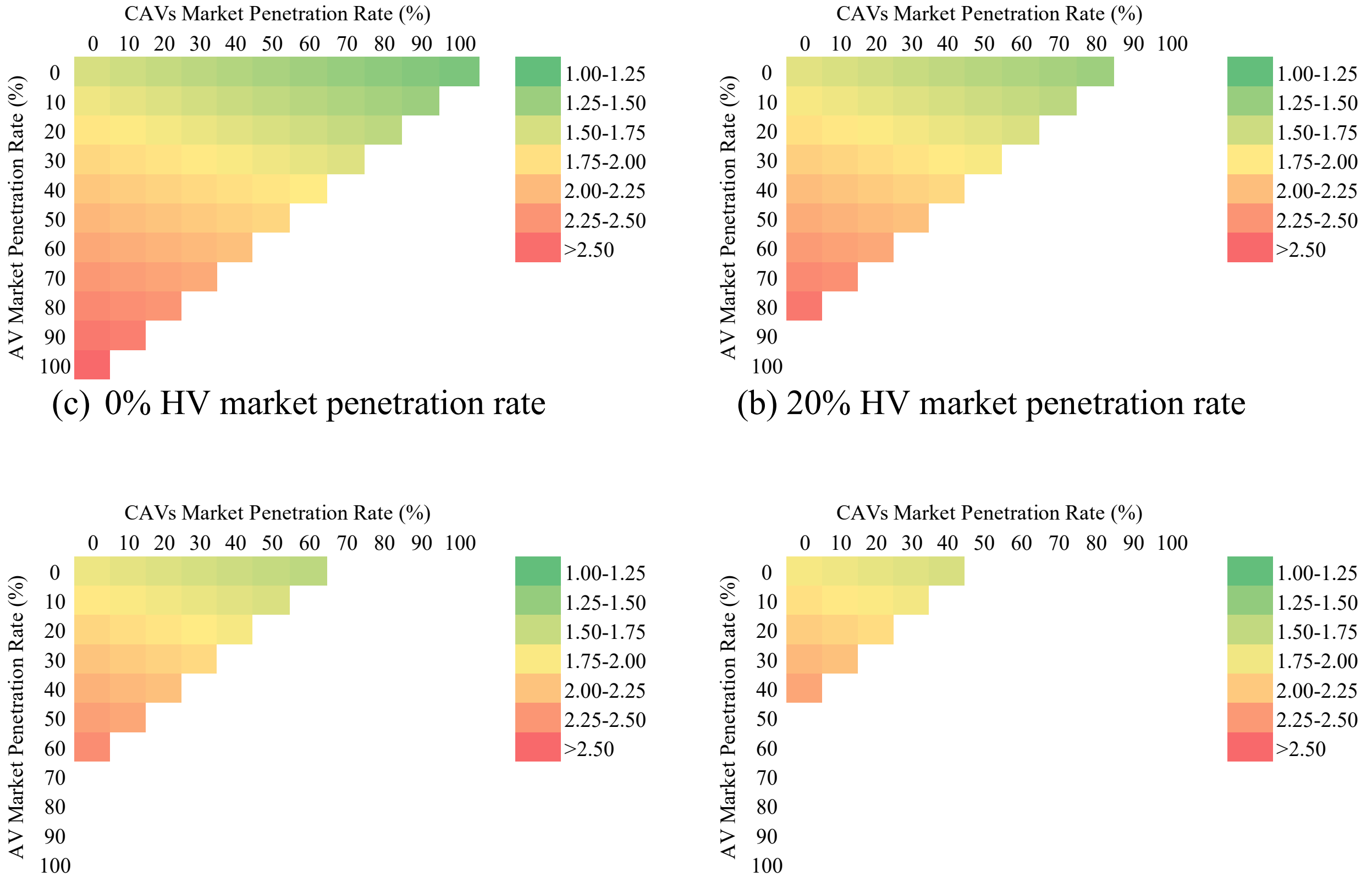

(c) 0% HV market penetration rate (b) 20% HV market penetration rate

(e) 40% HV market penetration rate (d) 60% HV market penetration rate

**FIGURE 11** Average intersection-level saturation headway (in seconds) for different CV, AV, and CAV market penetration rates. [Note: $CV = 1 - HV - AV - CAV$]

FIGURE 12 shows the the aggregated average saturation headway in the intersection. The saturation headway of human drivers was equal to two seconds, which was achieved by the calibration process. Increasing the CV market penetration rate decreased the average saturation headway to 1.5 seconds representing a 25% reduction, while increasing the AV market penetration rate increased it up to 2.6 seconds representing a 30% increase. CAVs moved more efficiently through the intersection with 1.2 seconds of saturation headway, indicating a 40% reduction compared to the base scenario.

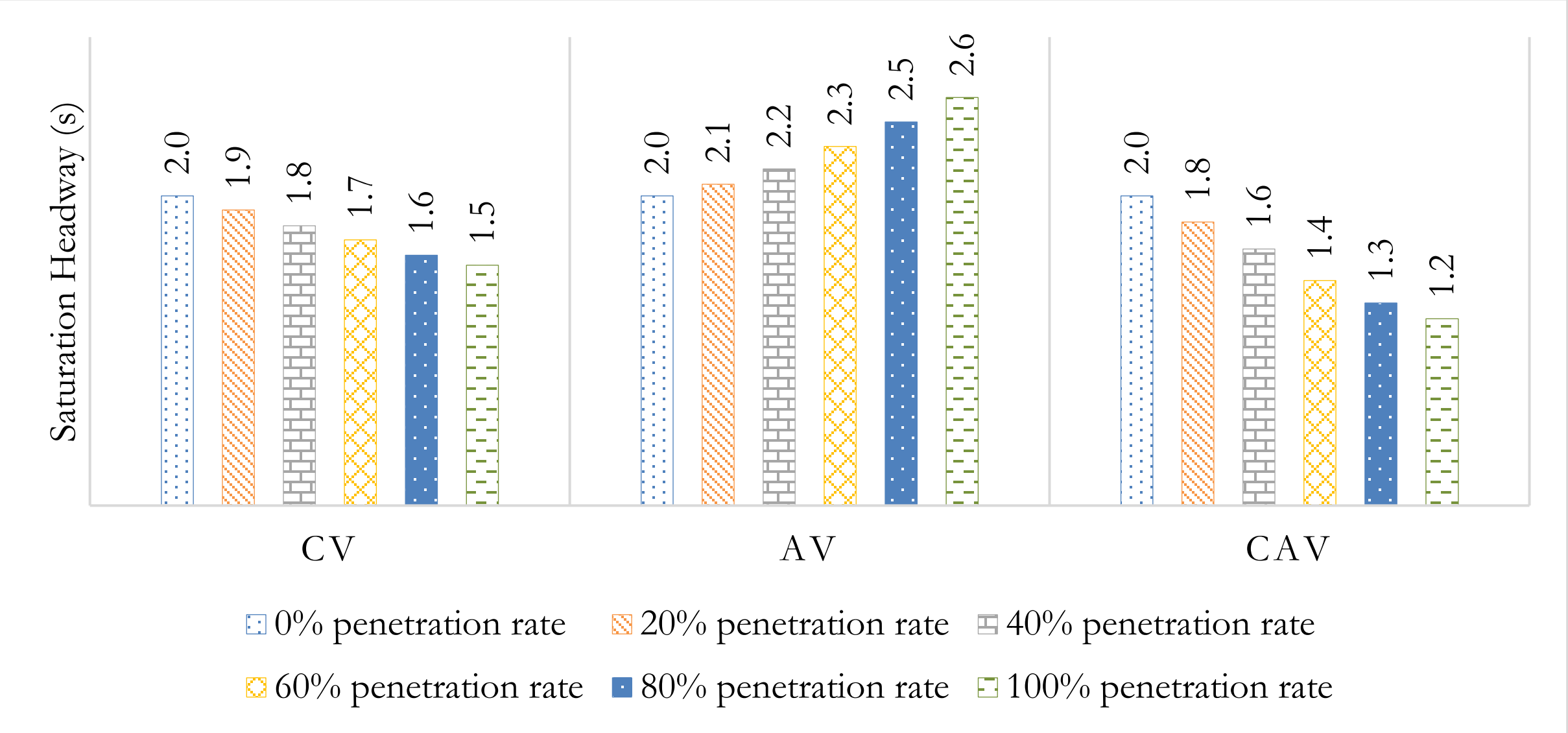

**FIGURE 12** Intersection-level saturation headway

## CONCLUSION

This study evaluates the potential effect of a mixed traffic stream of human-driven, automated, and connected vehicles on saturation headway and capacity at signalized intersections. Previous studies mainly focused on the operation of connected or automated vehicles on freeway facilities. However, the behavior of CVs, AVs, and CAVs in signalized intersections can be significantly affected by the information received on future signal timing plans. Four vehicle types are considered as (I) human-driven vehicles, (II) connected vehicles, (III) automated vehicles, and (IV) connected and automated vehicles. Vissim is used as a testbed to simulate the movement of vehicles with different driving behaviors and study their potential effects on mobility when they interact with each other and traffic signal controllers under various market penetration rates.

The results show that CAVs provide the most efficient mobility. CVs also improve mobility due to receiving information about the future signal timing plans. Both CVs and CAVs can adjust their speed upstream of the intersection to arrive during the green traffic light. In contrast with CVs and CAVs, AVs drive more cautiously and yield longer saturation headways and delays. Results also show that the highest increase (80%) and decrease (20%) in lane group capacity are observed respectively in a traffic stream of 100% CAVs and 100% AVs.

This study determines saturation headway and capacity adjustment factors for different lane groups under various CV, AV, and CAV market penetration rates. These values could be used to calculate the saturation flow rate and capacity of various lane groups in the presence of vehicles with full automation or connectivity cabilities. Both saturation headway and capacity adjustment factors are applicable to fixed-time and actuated controlled intersections. The results of this research can be used by transportation agencies to predict the impacts of connectivity and automation on local and regional traffic. This information is essential for long-range transportation plans that need to consider how to better prepare for the future and prevent or mitigate any adverse effects from emerging vehicle technologies. This study implements assumptions and changes, informed by the literature, in certain parameters of Vissim's car-following and lane-changing models, which were originally designed to represent human driving behavior. Further studies are required to replace the car following and lane changing logics of existing simulation packages with logics specifically designed for CVs, AVs, and CAVs. In addition, this study assumes conservative behavior for AVs. However, it is possible that manufactures will design AVs to operate less conservativly as the market penetration increases and the technology matures. Future studies can explore the potential changes in AV behavior over time and its implications on the capacity of signalized intersections.

## ACKNOWLEDGMENTS

This paper presents results from a larger project funded by the North Carolina Department of Transportation (NCDOT) titled "Impacts of Autonomous Vehicle Technology on Transportation Systems." The contents of this article reflect the views of the authors and do not reflect the official views or policies of NCDOT.

## AUTHOR CONTRIBUTION STATEMENT

The authors confirm contribution to the paper as follows: Ali Hajbabaie: Conceptualization, Data curation, Formal analysis, Investigation, Methodology, Project administration, Resources, Software, Supervision, Validation, Visualization, Writing - review & editing. Mehrdad Tajalli: Data curation, Formal analysis, Methodology, Validation, Visualization, Writing - original draft. Eleni Bardaka: Investigation, Project administration, Resources, Supervision, visualization, Writing - review & editing.